\documentclass{amsart}

\usepackage{amsmath, amscd, amssymb}
\usepackage[all]{xy}
\usepackage{color}

\newcommand{\half}{\frac{1}{2}}
\newcommand{\cal}{\mathcal}
\newcommand{\ba}{{\bf a}} \newcommand{\bb}{{\bf b}}
\newcommand{\pd}{\partial}

\newcommand{\bC}{{\Bbb C}}

\newcommand{\bQ}{{\Bbb Q}}
\newcommand{\bR}{{\Bbb R}}

\newcommand{\bZ}{{\Bbb Z}}

\newcommand{\cF}{{\cal F}}

\newcommand{\cO}{{\cal O}}
\newcommand{\cP}{{\cal P}}
\newcommand{\cX}{{\cal X}}

\newcommand{\cZ}{\mathcal{Z}}

\newcommand{\ma}{\mathfrak{a}} \newcommand{\mB}{\mathfrak{B}}

\newcommand{\vac}{|0\rangle}

 \DeclareMathOperator{\ch}{ch}

\DeclareMathOperator{\End}{End} 

\DeclareMathOperator{\Span}{span} \DeclareMathOperator{\id}{id}
 
 \DeclareMathOperator{\tr}{Tr}

\newcommand{\be}{\begin{equation}}
\newcommand{\ee}{\end{equation}}
\newcommand{\bea}{\begin{eqnarray}}
\newcommand{\ben}{\begin{eqnarray*}}
\newcommand{\een}{\end{eqnarray*}}
\newcommand{\eea}{\end{eqnarray}}
\newcommand{\bet}{\begin{equation}
\begin{split}}
\newcommand{\eet}{\end{split}
\end{equation}}

\newtheorem{theorem}{Theorem}[section]

\theoremstyle{remark}

\theoremstyle{definition}

\newtheorem{rmk}{Remark}

\definecolor{light-pink}{rgb}{1,.90,.90}
\definecolor{light-green}{rgb}{.95,1,.95}

\definecolor{yellow}{rgb}{1,1,0}
\definecolor{orange}{rgb}{1,.7,0}
\definecolor{red}{rgb}{1,0,0}
\definecolor{white}{rgb}{1,1,1}

\definecolor{A}{rgb}{.75,1,.75}

\newcommand{\corr}[1]{\langle {#1} \rangle}

\begin{document}

\title[Quasimodularity of Intersection Numbers on the Hilbert Schemes]
{On Quasimodularity of Some Equivariant Intersection Numbers on the Hilbert Schemes}

\author{Jian Zhou}
\address{Department of Mathematical Sciences\\
Tsinghua University\\ Beijing, China}
\email{jzhou@math.tsinghua.edu.cn}

\begin{abstract}
We  observe that certain equivariant intersection numbers
of Chern characters of tautological sheaves on Hilbert schemes for suitable
circle actions can be computed using the Bloch-Okounkov formula,
hence they are related to Gromov-Witten invariants of
elliptic curves and its operator formalism in terms of operators on the Fock space.
\end{abstract}

\maketitle

\section{Introduction}

Okounkov \cite{Oko} proposed several intriguing conjectures relating
intersection numbers of Chern classes on Hilbert schemes of points
to multiple q-zeta values.
For earlier results that lead to this conjecture, see \cite{Car1, Car2, Car-Oko}.
For recent results on this conjecture,
see e.g. \cite{Qin-Yu}.
For related works,
see e.g. \cite{Wan-Zho1, Wan-Zho2}.

In this paper
we will compute some $S^1$-equivariant intersection numbers on Hilbert schemes
of the affine plane,
and show that their suitably normalized generating series
are quasimodular forms.
We will also explain how to reduce their computations
to the Okounkov-Bloch formula \cite{BO},
hence a connection to the Gromov-Witten theory of elliptic curves \cite{Oko-Pan}
can be observed.
We will also compute the torus-equviariant intersection numbers and show that they 
can be reduced to the deformed Bloch-Okounkov $n$-point functions defined by Cheng and Wang \cite{Che-Wan}.

For the $S^1$-equivariant case, our idea is simply to combine the results Li-Qin-Wang \cite{LQW}
with equivariant localization and the theory of quasimodular forms \cite{Kan-Zag, BO, Zag}.
We reinterpret what has been computed in \cite[Theorem 4.5]{LQW} as
the trace of an operator on the   Fock space
as the generating series of equivariant intersection numbers
of equivariant Chern characters on the Hilbert schemes of points.
By localization,
equivariant intersection numbers on Hilbert schemes of points on $\bC^2$
are given by summations over fixed points
of torus actions on the Hilbert schemes,
which are indexed by partitions of integers.
Hence generating series of equivariant intersection numbers on $\bC^{[n]}$
are given by summations over partitions.
This is why they are related to the Bloch-Okounkov formula and quasimodularity.

Quasimodularity was suggested by Dijikgraaf \cite{Dij}
in mirror symmetry of elliptic curves.
The mathematical theory was formulated and developed
by Kaneko and Zagier \cite{Kan-Zag}.
Bloch and Okounkov \cite{BO} established the quasimodularity of summations of over partitions
for a large class of functions on partition function called shifted symmetric polynomials and derived an explicit formula.
Zagier \cite{Zag} developed a new approach to the results of Bloch and Okounkov.
Our note is a combination of some ideas from \cite{Car2, LQW}
with some ideas from \cite{BO, Zag}.
In an appendix we will present some new way to derive Zagier's Theorem 1 in \cite{Zag}.
Since the space of space of quasimodular forms
is contained in the space of multiple $q$-zeta values \cite{Oko},
we are dealing with a special case of Okounkov's conjecture.

For the torus-equivariant case,
the idea is similar. Equivariant localization reduces the computation
to a summation over partitions which we identify with deformed $n$-point 
functions defined by Cheng and Wang \cite{Che-Wan}
using a vertex operator realization of the Macdonald operator \cite{AMOS, Gar-Hai}.
Unfortunately closed formulas are available  in this case only for $n=1$ and $2$,
so a discussion of quasimodularity has to be left to further investigations at present.

We arrange the rest of this note as follows.
In \S 2 after introducing some notations on partitions,
we recall the work of Bloch and Okounkov \cite{BO}
on summation over partitions  and quasi-modular forms (see also Zagier \cite{Zag}).
We also make use of \cite[Lemma 3.1]{LQW} and reinterpret \cite[Theorem 4.5]{LQW}.
In \S 3 we recall the deformed $n$-point function and deformed
Bloch-Okounkov formula for $1$- and $2$-point functions \cite{Che-Wan}.
In \S 4 we explain how to reduce the computations of some equivariant
 intersection numbers of equivariant Chern characters
 of tautological sheaves on $\bC^2$ to Bloch-Okounkov formula or deformed Bloch-Okounkov characters.
In the Appendix we present a natural formalism that yields a formula that generalizes
both Jacobi's triple product identity and Zagier's
recursion relations for $q$-brackets of shifted symmetric polynomials.

\section{Partitions, Bloch-Okounkov Formula, and Quasimodular Forms}

In this section we first recall some notations on partitions
and recall a formula of Bloch-Okounkov \cite{BO}.
We also recall their work on quasimodularity.

\subsection{Notations on partitions}

A partition is a nonincreasing sequence $\lambda$ of nonnegative integers
$\lambda_1 \geq \lambda_2 \geq \cdots$,
with only finitely many zero terms.
One writes
$|\lambda|:=\sum_{i=1}^\infty \lambda_i$ and call it the weight of $\lambda$.
The length of the partition is defined to be the number of nonzero
terms in $\lambda$ and is denoted by $l(\mu)$.
The partition with length equal to $0$ is called the empty partition
and is denoted by $\emptyset$.

A partition $\lambda$ can be graphically represented by its Young diagram.
First one can assign the following numbers to  a box $s =(i,j)\in \lambda$:
\begin{align*}
a_{\lambda}(s) & = \lambda_i - j, & a_{\lambda}'(s) & = j-1,\\
l_{\lambda}(s) & = \lambda^t_j - i, & l_{\lambda}'(s) & = i -1.
\end{align*}
where $s$ is located at the $i$-row and the $j$-th column.
Note that
\begin{align*}
a_{\lambda^t}(s^t) & = l_{\lambda}(s), & a_{\lambda^t}'(s^t) = l_{\lambda}'(s),
\end{align*}
where $\lambda^t$ is obtained from $\lambda$ by switching the roles of rows and columns,
and $s^t$ is the box in $\lambda^t$ that corresponds to $s$ in $\lambda$.
The content of a box $s \in \lambda$ is defined by
\be
c(s):=j-i = a'_\lambda(s) - l'_\lambda(s).
\ee

\subsection{Some functions on partitions}

An  advantage of the graphical representation
by Young diagram is that one can naturally define
some functions on the set of partitions.

Recall the following formula \cite[Lemma 3.1]{LQW}:
\be \label{eqn:LQW}
\sum_{s\in \lambda} e^{zc_\lambda(s)}
= \frac{1}{\varsigma(z)}
\biggl( \sum_{i=1}^\infty e^{z(\lambda_i-i+1/2)}
- \frac{1}{\varsigma(z)} \biggr),
\ee
where  $\varsigma(z) = e^{z/2} - e^{-z/2}$.
The proof of the formula in \cite{LQW} is nice and simple.
For the convenience of the reader we recall their proof here.
The contents of $\lambda$ are:
$$-i + 1, -i + 2, \dots, -i + \lambda_i,\;\;\;
 i = 1, \dots, l(\lambda),$$
therefore,
\ben
\sum_{s\in \lambda} e^{zc_\lambda(s)}
& = & \sum_{i=1}^{l(\lambda)} (z^{-i + 1} +z^{-i + 2}+ \dots
+z^{-i + \lambda_i} ) \\
& = & \sum_{i=1}^{l(\lambda)} \frac{e^{z(\lambda_i-i+1)} - e^{z(-i+1)}}{e^z-1}\\
& = & \sum_{i=1}^{l(\lambda)} \frac{e^{z(\lambda_i-i+1/2)} - e^{z(-i+1/2)}}{e^{z/2}-e^{-z/2}} \\
& = & \frac{1}{\varsigma(z)}
\sum_{i=1}^\infty (e^{z(\lambda_i-i+1/2)} - e^{z(-i+1/2)}) \\
& = & \frac{1}{\varsigma(z)}
\biggl(\sum_{i=1}^\infty e^{z(\lambda_i-i+1/2)} - \sum_{i=1}^\infty e^{z(-i+1/2)}) \\
& = & \frac{1}{\varsigma(z)}
\biggl( \sum_{i=1}^\infty e^{z(\lambda_i-i+1/2)}
- \frac{1}{\varsigma(z)} \biggr).
\een
Here in the last equality the following identity has been used:
\be
\varsigma(z) = \sum_{i=1}^\infty e^{z(-i+1/2)}.
\ee
We will write the left-hand side of \eqref{eqn:LQW} as
$\ch_z(\lambda)$,
and call it the {\em Chern character} of $\lambda$
for reasons to be manifest in the next section.
It defines a function on the set $\cP$ of partitions.
For the empty partition $\emptyset$,
our convention is that $\ch_z(\emptyset) =1$.
We also write
\be
\ch_k(\lambda) = \frac{1}{k!} \sum_{c\in \mu} c_\lambda(s)^k.
\ee
From our convention for $\ch_z(\emptyset)$,
we have
\be
\ch_k(\emptyset) = \delta_{k,0}.
\ee

Another advantage of using the Young diagram is that one can introduce
the Frobenius notation $(a_1, \dots, a_r|b_1, \dots, b_r)$ graphically,
where $a_i$ (resp. $b_i)$
are the numbers of the cells to the right of (resp. below) the $i$-th cells on the diagonal.
An alternative way to understand these numbers is to use the Dirac sea.
Consider the set $X_\lambda:=\{\lambda_i-i+1/2|i=1,2,\dots\}$.
Denote by $(\bZ+\half)_\pm$ the positive (resp. negative) half integers,
then one has (see e.g. \cite[(21)]{Zag}):
\ben
&& X_\lambda \cap (\bZ+\half)_+ = C_\lambda^+:= \{a_r+\half, \dots, a_1+ \half\}, \\
&& (\bZ+\half) -X_\lambda = C_\lambda^-:=\{-b_1-\half, \dots, -b_r-\half\}.
\een
One can then define a sequence $P_k(\lambda)$ of functions of partitions for $k \geq 0$ as follows:
\be
P_k(\lambda) =  \sum_{j=1}^r ((a_j+\half)^k - (-b_j-\half)^k).
\ee
It is easy to see that
\be
P_0(\lambda) = 0.
\ee

One can relate these functions to the Chern characters of partitions.
From \eqref{eqn:LQW} we have
\be
\varsigma(z)\ch_z(\lambda) = \sum_{k=1}^\infty \frac{z^k}{k!} P_k(\lambda).
\ee
Indeed, this follows from the following straightforward computations:
\ben
\varsigma(z)\ch_z(\lambda) & = & \varsigma(z)\sum_{s\in \lambda} e^{zc_\lambda(s)}
=   \sum_{i=1}^\infty e^{z(\lambda_i-i+1/2)}
- \frac{1}{\varsigma(z)} \\
& = & \sum_{\lambda_i-i+1/2 \in C_\lambda^+} e^{z(\lambda_i-i+1/2)}
- \sum_{\lambda_i-i+1/2 \in C_\lambda^-} e^{z(\lambda_i-i+1/2)} \\
& = & \sum_{j=1}^r (e^{z(a_j+\half)} - e^{-z(b_j+\half)}  )\\
& = & \sum_{j=1}^r \sum_{k=0}^\infty \frac{z^k}{k!} ((a_j+\half)^k - (-b_j-\half)^k) \\
& = & \sum_{k=0}^\infty \frac{z^k}{k!} P_k(\lambda).
\een
By comparing the coefficients of $z^k$,
one then gets
\be
P_k(\lambda) = k! \sum_{2l+m=k} \frac{2^{1-2l}}{(2l)!} \ch_m(\lambda).
\ee
One also has
\be
\ch_z(\lambda) =
\frac{1}{e^{z/2}-e^{-z/2}} \sum_{k=1}^\infty \frac{z^k}{k!} P_k(\lambda).
\ee
Write
\be
\frac{z/2}{e^{z/2}-e^{-z/2}} = \sum_{n=0}^\infty \beta_n z^n
\ee
as in \cite{Zag}, where
\be
\beta_n = \frac{B_n}{n!}(\frac{1}{2^n} - \frac{1}{2}) = \half \frac{B_n(1/2)}{n!}.
\ee
Then one can get
\be
\ch_n(\lambda) = \sum_{k+l=n+1} 2\beta_k \frac{P_l(\lambda)}{l!}.
\ee

Zagier \cite{Zag} introduced another sequence of functions $Q_n$ of partitions:
\be
\sum_{i=1}^\infty e^{z(\lambda_i-i+1/2)} = \sum_{k=0}^\infty Q_k(\lambda) z^{k-1}.
\ee
Since we have
\ben
\sum_{i=1}^\infty e^{z(\lambda_i-i+1/2)}
& = & \sum_{k=0}^\infty \frac{z^k}{k!} P_k(\lambda) + \frac{1}{e^{z/2}-e^{-z/2}} \\
& = & \sum_{k=0}^\infty \frac{z^k}{k!} P_k(\lambda)
+ 2\sum_{k=0}^\infty \beta_k z^{k-1},
\een
the following equality holds:
\be
Q_k(\lambda) = \frac{P_{k+1}(\lambda)}{(k+1)!} + 2\beta_k
= \frac{P_{k+1}(\lambda)}{(k+1)!} + \frac{B_k(1/2)}{k!}.
\ee

\subsection{Bloch-Okounkov formula}

Let $f(\lambda)$ be a function on partitions,
Bloch and Okounkov \cite{BO} define
\be
\corr{f}_q :
= \frac{\sum_{\lambda} f(\lambda)q^{|\lambda|}}{\sum_{\lambda} q^{|\lambda|}}.
\ee
For $n=1,2,3, \dots,$,
they also define
\be
F(t_1, \dots, t_n; q):
= \biggl\langle \prod_{k=1}^n \biggl( \sum_{i=1}^\infty t_k^{\lambda_i-i+\frac{1}{2}} \biggr)
\biggr\rangle_q,
\ee
and prove the following famous formula:
\be \label{eqn:BO}
F(t_1, \dots, t_n)
= \sum_{\sigma \in S_n}
\frac{\det \biggl(
\frac{\theta^{(j-i+1)}(t_{\sigma(1)} \cdots t_{\sigma(n-j)}) } {(j-i+1)!} \biggr)}
{\theta(t_{\sigma(1)}) \cdots \theta(t_{\sigma(1)}t_{\sigma(2)}) \cdots
\theta(t_{\sigma(1)} \cdots t_{\sigma(n)})}.
\ee
Here $\theta(x)$ is the Jacobi theta function defined by
\be
\theta(t):= \eta(q)^{-3} \sum_{n\in \bZ}
(-1)^nx^{n+ 1/2}q^{\frac{(n+ 1/2)^2}{2}}
 = (q)_\infty^{-2}(t^{1/2} - t^{-1/2})(qt)_\infty(q/t)_\infty,
\ee
and $\theta^{(p)}(t) = (t\frac{d}{dt})^p\theta(t)$.
See \cite[\S 5]{Oko-Pan} for the relationship of this formula
to the stationary Gromov-Witten invariants of elliptic curves.

Combining formula \eqref{eqn:LQW} with \eqref{eqn:BO},
one can easily compute the $n$-point function of
Chern characters of partitions:
\be
\langle \ch_{z_1} \cdots \ch_{z_n} \rangle_q
=\biggl\langle \prod_{j=1}^n \frac{1}{\varsigma(z_j)}
\biggl( \sum_{i=1}^\infty e^{z_j(\lambda_i-i+1/2)}
- \frac{1}{\varsigma(z_j)} \biggr) \biggr\rangle_q
\ee
by expanding the right-hand side:
\be \label{eqn:Chern}
\begin{split}
& \langle \ch_{z_1} \cdots \ch_{z_n} \rangle_q \\
= & \frac{1}{\prod_{j=1}^n \varsigma(z_j)} \sum_{k=0}^n
\sum_{1 \leq i_1 < \cdots < i_k \leq n}
F(e^{z_{i_1}}, \dots, e^{z_{i_k}})
\frac{(-1)^{n-k}\varsigma(z_{i_1}) \cdots \varsigma(z_{i_k})}
{\prod_{j=1}^n \varsigma(z_j)}.
\end{split}
\ee
See \cite[Theorem 4.5]{LQW} for the interpretation of this formula
in terms of the trace of an operator on the Fock space.

\subsection{Quasimodular forms}

Quasimodularity was suggested by Dijkgraaf \cite{Dij} in the context of mirror symmetry
of elliptic curves.
The mathematical theory of quasimodular forms was developed by Kaneko and Zagier \cite{Kan-Zag}.
The ring of quasimodular forms is $\bQ[E_2, E_4, E_6]$.
The work of Bloch and Okounkov \cite{BO} provides many quasimodular forms.
They have shown that for any polynomial $f$ in $P_1, P_2, \dots$,
$\corr{f}_q$ is a quasimodular form.
For a different proof of this result,
see Zagier \cite{Zag}.
As a consequence,
$\corr{\ch_{k_1} \cdots \ch_{k_n}}_q$ is a quasimodular form.

\section{Deformed Bloch-Okounkov Formula}

In this section we recall the work on deformed Bloch-Okounkov formula due to Cheng and Wang \cite{Che-Wan}.
We follow their notations and presentation closely with only minor modifications.

\subsection{Deformed vertex operator}
Consider the Heisenberg algebra generated by $I$ and $\ma_n$, $n \in \bZ$,
with the commutation relations:
$$[\ma_m, \ma_n] = \kappa m \delta_{m,-n}I,$$
where $\kappa$ is a Planck constant.
The bosonic Fock space $B$ has a basis
$\ma_{-\lambda}:= \ma_{-\lambda_1}\ma_{-\lambda_2} \cdots \vac$,
where $\lambda = (\lambda_1, \lambda_2, \dots)$ runs over all partitions.
After identifying it with the ring of symmetric function $\Lambda$
by identifying $\ma_{-\lambda}$ with the power-sum symmetric functions $p_\lambda$,
one has ($n \geq 1$):
\begin{align}
\ma_{-n} & = p_n \cdot, & \ma_0 & = 0, & \ma_n & = \kappa n \frac{\pd}{\pd p_n}, & I & =1.
\end{align}
Denote by $\Lambda_{q,t}$ the ring of symmetric functions with coefficients in $\bQ(q, t)$.
Denote by $L_0$ the usual energy operator on $\Lambda_{q,t}$:
$L_0g = n g$ if $g$ is a symmetric function of degree $n$.
For $f\in \End(\Lambda_{q,t}$,
the trace of $A$ is defined by:
\be
\tr_vA:= \tr(v^{L_0}A)|_{\Lambda_{q,t}}.
\ee
As in \cite{BO}, for a function $f$ on the set $\cP$ of partitions,
define the $v$-bracket of $f$ by
\be
\corr{f}_v: = (v)_\infty \sum_{\lambda \in \cP} f(\lambda) v^{|\lambda|}.
\ee
Cheng and Wang introduced the following deformed vertex operator:
\be
V(z; q_1, t_1, q_2, t_2) = \exp \biggl( \sum_{k\geq 1}
(q_1^k - q_2^k)a_{-k} \frac{z^k}{k} \biggr) \cdot
\exp \biggl(- \sum_{k\geq 1} (t_1^k - t_2^k) \alpha_k\frac{z^{-k}}{k}\biggr).
\ee
They proved the following formula \cite[Theorem 16]{Che-Wan}:
\be
\begin{split}
\biggl\langle \prod_{i=1}^n V(z_i; s_i,t_i, u_i, w_i) \biggr\rangle_v
= & \prod_{1 \leq i < j \leq n}
\biggl[ \frac{(1-t_i w_jz_i^{-1}z_j)(1-s_iu_jz_i^{-1}z_j)}
{(1-t_i u_jz_i^{-1}z_j)(1-s_iw_jz_i^{-1}z_j)} \biggr]^\kappa \\
& \cdot \prod_{i,j=1}^n \biggl[ \frac{(t_i w_jz_i^{-1}z_jv)_\infty(s_iu_jz_i^{-1}z_jv)_\infty}
{(t_i u_jz_i^{-1}z_jv)_\infty(s_iw_jz_i^{-1}z_jv)_\infty} \biggr]^\kappa,
\end{split}
\ee
where $(a)_\infty : = \prod_{i=0}^\infty (1-av^i)$.
For example,
when $n=1$,
\be \label{eqn:1-point-V}
\langle V(z_1; s_1,t_1, u_1, w_1) \rangle_v
= \biggl[ \frac{(t_1 w_1v)_\infty(s_1u_1v)_\infty}{(t_1 u_1v)_\infty(s_1w_1v)_\infty} \biggr]^\kappa.
\ee

Write
\be
V(z; q_1, t_1, q_2, t_2) = \sum_{m\in \bZ}
V_m(q_1, q_2, t_1, t_2)z^m.
\ee
The operator $V_0$ is called the zero mode of $V$.
By \eqref{eqn:1-point-V} one can get:
\be \label{eqn:1-point}
\langle V_0(s_1,t_1, u_1, w_1) \rangle_v
= \biggl[ \frac{(t_1 w_1v)_\infty(s_1u_1v)_\infty}{(t_1 u_1v)_\infty(s_1w_1v)_\infty} \biggr]^\kappa.
\ee
This is \cite[Theorem 13]{Che-Wan} proved by a different method.

\subsection{Deformed Bloch-Okounkov formula}

As pointed out by Cheng and Wang \cite[Remark 12]{Che-Wan},
when $\kappa = 1$, $q_2 = t_2 = 1$, and write $q = q_1$ and $t = t_1$,
the operator $V_0$ provides a vertex operator realization for
the  Macdonald operator $\hat{\mB}_{q,t}$:
\be
\hat{\mB}_{q,t} =
\frac{1}{(1 - q)(1 - t)} \cdot V_0(q, 1, t, 1).
\ee
This formula appears in a different form in the study of Macdonald
polynomials by Garsia and Haiman \cite[(73)]{Gar-Hai}.
It also appears in \cite[(32)]{AMOS}.
This operator has Macdonald functions $P_\lambda(x; q, t)$ as eigenfunctions:
\bea
&& \hat{\mB}_{q,t} P_\lambda(x;q,t) = \hat{B}_\lambda(q,t) \cdot P_\lambda(x;q,t), \\
&& \hat{B}_\lambda(q,t) : = \frac{1}{1-q} \sum_{i\geq 1} t^{i-1} q^{\lambda_i}.
\eea
There is a related operator:
\bea
&& \mB_{q,t} P_\lambda(x;q,t) = B_\lambda(q,t) \cdot P_\lambda(x;q,t), \\
&& B_\lambda(q,t) : = \sum_{\Box \in \lambda} q^{a'(\Box)} t^{l'(\Box)}. \label{eqn:Blambda}
\eea

The $n$-point (correlation) functions are defined to be
\bea
&& F(q_1, t_1; \dots ; q_n, t_n) := \tr_v(\mB_{q_1,t_1}, \dots,  \mB_{q_n,t_n}), \\
&& \hat{F}(q_1, t_1; \dots ; q_n, t_n) := \tr_v(\hat{\mB}_{q_1,t_1}, \dots,  \hat{\mB}_{q_n,t_n}).
\eea
By \cite[Lemma 1]{Che-Wan},
\be
B_\lambda(q,t) = \hat{B}_\emptyset(q,t) - \hat{B}_\lambda(q,t)
= \frac{1}{(1-q)(1-t)},
\ee
so one convert between the computations of $F$ and that of $\hat{F}$.
By \cite[Lemma 2]{Che-Wan},
\be
F(q_1, t_1; \dots , q_n, t_n) = (v)_\infty^{-1} \biggl\langle \prod_{k=1}^n B_\lambda(q_k, t_k)\biggr\rangle_v.
\ee
The one-point function and the two-point function have been computed by Cheng and Wang
\cite[Theorem 5, Theorem 9]{Che-Wan}:
\bea
&& \hat{F}(q,t) = \frac{(vqt)_\infty}{(q)_\infty(t)_\infty}, \\
&& \hat{F}(q_1,t_1; q_2, t_2) = \frac{1}{(1-q_1)(1-q_2)(1-t_1t_2)}
\cdot \frac{(vq_1q_2t_1t_2)_\infty}{(vt_1t_2)_\infty(vq_qq_2)_\infty} \\
&& \cdot \biggl[ \frac{q_1q_2t_1t_2-1}{(1-q_1t_1)(1-q_2t_2)}
+ \frac{1}{1-q_1t_1}\Phi(v, q_1t_1, vq_1q_2; vq_1, vq_1q_2t_1t_2; v; t_2) \nonumber \\
& & +  \frac{1}{1-q_1t_1}\Phi(v, q_2t_2, vq_1q_2; vq_2, vq_1q_2t_1t_2; v; t_1) \biggr], \nonumber
\eea
where $\Phi$ is $v$-basic hypergeometric series:
\be
\Phi(a_1, a_2, a_3; b_1, b_2; z)
: = \sum_{m \geq 0} \frac{(a_1)_m(a_2)_m(a_3)_m}{(b_1)_m(b_2)_m} \frac{z^m}{(v)_m},
\ee
 and $(a)_0=1$, $(a)_m= \prod_{i=0}^{m-1} (1-av^i)$.

\section{Equivariant Intersection Numbers on $\bC^{[n]}$ and Quasimodularity}

See \cite{Nak} for references on equivariant indices
on Hilbert schemes of points on $\bC^2$.
We have followed the notations in \cite{Wan-Zho1, Wan-Zho2}.

\subsection{Localizations on Hilbert schemes of the affine plane}

By a theorem of Fogarty \cite{For} the Hilbert scheme $(\bC^2)^{[n]}$
is a nonsingular variety of dimension $2n$.
The torus action on $\bC^2$ given by
$$(t_1, t_2) \cdot x = t_1 x, \;\; (t_1, t_2) \cdot y = t_2 y$$
on linear coordinates induces an action on $(\bC^2)^{[n]}$.
The fixed points are isolated and parameterized by partitions
$\lambda=(\lambda_1, \dots, \lambda_l)$ of weight $n$.
They correspond to ideals
$$I_{\lambda} = \langle y^{\lambda_1}, xy^{\lambda_2}, \dots, x^{l-1}y^{\lambda_l}, x^l\rangle.$$
The weight decomposition of the tangent bundle of $T(\bC^2)^{[n]}$ at a
fixed point $\lambda$ is given by \cite{Nak}:
\begin{eqnarray} \label{eqn:Weights}
&& \sum_{(i, j)\in \lambda} (t_1^{(\lambda^t_j - i)} t_2^{-(\lambda_i-j+1)}
+ t_1^{-(\lambda^t_j-i+1)}t_2^{(\lambda_i-j)}) \\
& = & \sum_{s \in \lambda} (t_1^{l(s)} t_2^{-(a(s)+1)}
+ t_1^{-(l(s)+1)}t_2^{a(s)}),
\end{eqnarray}
where $(t_1, t_2)\in \bC^*\times \bC^*$.
It follows that the equivariant Euler class is given at $I_\lambda$ by:
\be
\begin{split}
e_T(T\bC^{[n]})|_{I_\lambda} & = ((\lambda^t_j - i)t_1 -(\lambda_i-j+1)t_2)\cdot
(-(\lambda^t_j-i+1)t_1+(\lambda_i-j)t_2) \\
& = \prod_{s \in \lambda} (l(s)t_1-(a(s)+1)t_2) (-(l(s)+1)t_1 +a(s)t_2).
\end{split}
\ee

\subsection{Tautological bundles}

Let $\cZ_n\subset X\times X^{[n]}$ be the universal family of subschemes parameterized by
$X^{[n]}$. Denote by $p_1: \cZ_n \to X$ and $\pi: \cZ_n \to X^{[n]}$
the projection onto the $X$ and $X^{[n]}$ respectively.
For any locally free sheaf $F$ on $X$ let $F^{[n]} = \pi_*(\cO_{\cZ_n}
\otimes p_1^*F)$.
With this notation we write $\xi_n=\xi_n^X=\cO_X^{[n]}$.
The tautological bundle $\xi_n$ on $(\bC^2)^{[n]}$
has its weight decomposition at a fixed point $I_\lambda$ given by \cite{Nak}:
\begin{eqnarray*}
&& \xi_n|_{I_{\lambda}} =  \sum_{(i, j) \in \lambda} t_1^{i-1}t_2^{j-1}
= \sum_{s \in \lambda} t_1^{l'(s)}t_2^{a'(s)}.
\end{eqnarray*}
So we have
\be
\begin{split}
\ch(\xi_n)_T|_{I_{\lambda}} & = \sum_{(i, j) \in \lambda} e^{(i-1)t_1+(j-1)t_2}  \\
& = \sum_{s \in \lambda} e^{l'(s)t_1+a'(s)t_2}.
\end{split}
\ee
In particular,
\be
\begin{split}
\ch_k(\xi_n)_T|_{I_{\lambda}} & =  \frac{1}{k!} \sum_{(i, j) \in \lambda} ((i-1)t_1+(j-1)t_2)^k  \\
& = \frac{1}{k!} \sum_{s \in \lambda} (l'(s)t_1+a'(s)t_2)^k.
\end{split}
\ee

For a vector $A =(a, b) \in \bZ^2$,
denote by $\cO_{\bC^2}^A$ the $T^2$-equivariant line bundle on $\bC^2$
with weight $A$.
Recall the universal family $\cZ_n$ lies in $\bC^2 \times (\bC^2)^{[n]}$,
and denote by $p_1: \cZ_n \to \bC^2$ the projection onto the first factor.
Let $\xi_n^A = \pi_*(\cO_{\cZ_n} \otimes p_1^*\cO_{\bC^2}^A)$.
Then one has:
\begin{eqnarray*}
&& \xi_n^A|_{I_{\lambda}} =  \sum_{(i, j) \in \lambda}
t_1^{i-1}t_2^{j-1}t_1^{a}t_2^{b} = \sum_{s \in \lambda}
t_1^{l'(s)}t_2^{a'(s)}t_1^{a}t_2^{b}.
\end{eqnarray*}
So we have:
\be
\begin{split}
\ch(\xi_n^A|_{I_{\lambda}})_T & =  \sum_{(i, j) \in \lambda} e^{(i-1+a)t_1+(j-1+b)t_2}  \\
& = \sum_{s \in \lambda} e^{(l'(s)+a)t_1+(a'(s)+b)t_2}.
\end{split}
\ee
In particular,
\be
\begin{split}
\ch_k(\xi_n^A|_{I_{\lambda}})_T & =  \frac{1}{k!} \sum_{(i, j) \in \lambda} ((i-1+a)t_1+(j-1+b)t_2)^k  \\
& = \frac{1}{k!} \sum_{s \in \lambda} ((l'(s)+a)t_1+(a'(s)+b)t_2)^k.
\end{split}
\ee

Hence by localization formula, we have
\begin{equation} \label{eqn:LefMain}
\begin{split}
& \sum_{n \geq 0} q^n
\int_{(\bC^2)^{[n]}_T} \ch_{k_1}(\xi_n^{A_1})_T \cdots \ch_{k_N}(\xi_n^{A_n})_T \cdot e_T(T\bC^{[n]})  \\
= & \sum_{\lambda} q^{|\lambda|}
\prod_{j=1}^N \frac{1}{k_j!} \sum_{s \in \lambda} ((l'(s)+a_j)t_1+(a'(s)+b_j)t_2)^{k_j}.
\end{split}
\end{equation}

Consider the circle subgroup $S^1 \to T$, $e^{it} \mapsto (e^{-it}, e^{it})$.
I.e., let $-t_1=t_2 = t$.
\ben
&& \sum_{n \geq 0} q^n
\int_{(\bC^2)^{[n]}_{S^1}} \ch_{k_1}(\xi_n^{A_1})_{S^1} \cdots
\ch_{k_N}(\xi_n^{A_n})_{S^1} \cdot e_{S^1}(T\bC^{[n]})  \\
& = & t^{k_1+\cdots + k_N} \sum_{\lambda} q^{|\lambda|}
\prod_{j=1}^N \frac{1}{k_j!} \sum_{s \in \lambda} (c(s)+(a_j-b_j))^{k_j} \\
& = & t^{k_1+\cdots + k_N} \sum_{\lambda} q^{|\lambda|}
\prod_{j=1}^N \sum_{l_j=0}^{k_j} \sum_{s \in \lambda} \frac{1}{l_j!} c(s)^{l_j} \cdot
\frac{1}{(k_j-l_j)!} (a_j-b_j)^{k_j-l_j}.
\een
The right-hand side is a linear combination of terms of the form:
\ben
t^{k_1+\cdots + k_N} \sum_{\lambda} q^{|\lambda|}
\prod_{j=1}^N \sum_{s \in \lambda} c(s)^{l_j}.
\een
Therefore,
by the results of last section,
$\frac{1}{\prod_{n=1}^\infty (1-q^n)} \sum_{\lambda} q^{|\lambda|}
\prod_{j=1}^N \sum_{s \in \lambda} c(s)^{l_j}$ is a quasimodular form.
Therefore, we have proved:

\begin{theorem}
For integers $k_1, \dots, k_N \geq 0$, $A_1, \dots, A_N \in \bZ^2$,
$$\frac{1}{t^{k_1+\cdots +k_N} \prod_{n=1}^\infty (1-q^n)}
\sum_{n \geq 0} q^n
\int_{(\bC^2)^{[n]}_{S^1}} \ch_{k_1}(\xi_n^{A_1})_{S^1} \cdots
\ch_{k_N}(\xi_n^{A_N})_{S^1} \cdot e_{S^1}(T\bC^{[n]})
$$
is a quasimodular form.
\end{theorem}

In a similar way, one has
\ben
&& \frac{1}{\prod_{n=1}^\infty (1-q^n)}
\sum_{n \geq 0} q^n
\int_{(\bC^2)^{[n]}_{S^1}} \ch(\xi_n^{A_1})_{S^1} \cdots
\ch(\xi_n^{A_N})_{S^1} \cdot e_{S^1}(T\bC^{[n]}) \\
& = & \frac{e^{\sum_{j=1}^N (b_j-a_j)t}}{\prod_{n=1}^\infty (1-q^n)} \sum_{\lambda} q^{|\lambda|}
\ch_t(\lambda) \cdots \ch_t(\lambda) \\
& = & e^{\sum_{j=1}^N (b_j-a_j)t} \cdot \corr{\ch_{t} \cdots \ch_{t}}_q.
\een
Hence it is possible to reduce to the Bloch-Okounkov formula by \eqref{eqn:Chern}.
In order to  compute the series in the above Theorem,
one can consider their generating series:
\ben
&& \frac{1}{\prod_{n=1}^\infty (1-q^n)}
\sum_{n \geq 0} q^n
\int_{(\bC^2)^{[n]}_{S^1}} \prod_{j=1}^N
\sum_{k_j\geq0} z_j^{k_j}\ch_{k_j}(\xi_n^{A_j})_{S^1} \cdot e_{S^1}(T\bC^{[n]}) \\
& = & \frac{e^{\sum_{j=1}^N (b_j-a_j)z_jt}}{\prod_{n=1}^\infty (1-q^n)} \sum_{\lambda} q^{|\lambda|}
\ch_{z_1t}(\lambda) \cdots \ch_{z_Nt}(\lambda) \\
& = & e^{\sum_{j=1}^N (b_j-a_j)z_jt} \cdot \corr{\ch_{z_1t} \cdots \ch_{z_Nt}}_q.
\een

For the torus-equivarinat intersection numbers, we consider:
\ben
&& \sum_{n \geq 0} q^n
\int_{(\bC^2)^{[n]}_{S^1}} \prod_{j=1}^N
\sum_{k_j\geq0} z_j^{k_j}\ch_{k_j}(\xi_n^{A_j})_{T} \cdot e_{T}(T\bC^{[n]}) \\
& = & \sum_{\lambda} q^{|\lambda|} \prod_{k=1}^N \sum_{(i, j) \in \lambda}
e^{z_k[(i-1+a_k)t_1+(j-1+b_k)t_2]} \\
& = & e^{\sum_{j=1}^N (a_jt_1+b_jt_2)z_j} \cdot 
\sum_{\lambda} q^{|\lambda|} \prod_{k=1}^N \sum_{s \in \lambda}
e^{z_k[l'_\lambda(s)t_1+a'_\lambda(s)t_2]}.
\een
By \eqref{eqn:Blambda},
its computation is reduced to the deformed $n$-point function
\be
\corr{\mB_\lambda(e^{z_1t_2}, e^{z_1t_1}) \cdots \mB_\lambda(e^{z_Nt_2}, e^{z_Nt_1})}_q.
\ee
This solves the problem of finding geometric interpretations of such $n$-point
functions posed by Cheng and Wang in the end of their paper \cite{Che-Wan}.

\begin{appendix}

\section{Spectral Flow, Dirac Flow, and Zagier Recursion Relations}

In this Appendix we interpret Zagier's approach in \cite{Zag}
in terms of spectral flow and Dirac flow,
and derive a formula that generalizes both the Jacobi triple product identity
and Zagier's recursion relations \cite[Theorem 1]{Zag}.

\subsection{The spectral flow}

For $\alpha \in \bR$,
note
\be
 D e^{2\pi i(n+\alpha)x} = (n+\alpha) \cdot e^{2\pi i (n+\alpha)x},
\ee
where $D:=\frac{1}{2\pi i}\frac{d}{dx}$ can be understood as the Dirac operator
on the line,
acting on the spaces of quasi-periodic functions:
\be
f(x+1) = e^{2\pi i \alpha} \cdot f(x).
\ee
The functions $e_{n+\alpha}(x) := e^{2\pi i (n+\alpha)x}$ are
the eigenfunctions in this space,
and one can use them to form the ``Dirac sea" as follows.

Denote by $(\bZ+\alpha)_\pm$ the nonnegative (resp. negative) numbers in the set
$\bZ + \alpha: =\{n+\alpha\;|\; n\in \bZ\}$.
It is clear that when $\alpha -\alpha' \in \bZ$,
then $\bZ+\alpha = \bZ+\alpha'$,
and so $(\bZ+\alpha)_\pm = (\bZ+\alpha')_\pm$.
It then suffices to consider $\alpha \in [0,1]$.
For a sequence $\ba = \{a_1, a_2, \dots\} \subset \bZ +\alpha$
with $a_1 > a_2 > \cdots$,
define
\begin{align}
\ba_+ & := \ba \cap (\bZ+\alpha)_+, &
\ba_-^c & := (\bZ+\alpha)_--\ba.
\end{align}
When both $\ba_+$ and $\ba^c_-$ are finite sets,
we say $\ba$ is {\em admissible}.
Denote by $\cX_\alpha$ the set of all admissible sequences in $\bZ+\alpha$.
To each $\ba \in \cX_\alpha$,
define
\be
e_{\ba} : = e_{a_1} \wedge e_{a_2} \wedge \cdots,
\ee
and define the (fermionic) Fock space $\cF_\alpha$ by
\be
\cF_\alpha:=\Span \{e_\ba\;|\; \ba \in \cX_\alpha\}.
\ee
Two special cases have been widely used in the mathematical literature:
$\cF_0$ (the Ramond-Ramond sector) and $\cF_{1/2}$ (the Neveu-Schwarz sector).
But the spectral flow in the physics literature \cite{SS, LVW} indicates that
it is interesting to relate these two sectors by continuously going through other sectors.
In our case,
this can be easily done by introducing a shift operator.
For $\beta \in \bR$,
define $U(\beta): \bZ + \alpha \to \bZ+\alpha + \beta$ simply by
$$n+\alpha \mapsto n+\alpha+\beta.$$
The family $\{U(\beta)\}_{\beta\in \bR}$ is called the {\em spectral flow}.
First of all, if is a flow because:
\begin{align}
U(0) & = \id, & U(\beta_1) U(\beta_2) & = U(\beta_1+\beta_2);
\end{align}
secondly,
it flows the eigenfunctions and the eigenvalues.
The operator $U(\beta)$ induces a shift operator, also denoted by $U(\beta)$, from
$\cX_\alpha$ to $\cX_{\alpha+\beta}$:
$$\ba \in \cX_\alpha \mapsto \ba + \beta \in \cX_{\alpha+\beta},$$
where $\ba+\beta=\{a_1+\beta, a_2+\beta, \dots\}$ when
$\ba =\{a_1, a_2, \dots\}$.
One therefore also has an induced operator, also denoted by $U(\beta)$,
from $\cF_\alpha$ to $\cF_{\alpha+\beta}$,
defined by $e_\ba\mapsto e_{U(\beta)e_\ba}$.

The following is the starting point of the approach of Zagier \cite{Zag}.
Define the charge of $e_\ba$ by $|\ba_+|-|\ba^c_-|$.
Inside $\cF_{1/2}$,
the charge zero $e_\ba$ are in one-to-one correspondence with the set of partitions.
Given a partition $\lambda=(\lambda_1, \lambda_2, \dots)$,
define
\be
\ba_\lambda = (\lambda_1-1+1/2, \lambda_2-2+1/2, \dots),
\ee
then $e_{\ba_\lambda}$ has charge zero in $\cF_{1/2}$,
and all $e_\ba$ in $\cF_{1/2}$ of charge zero has this form.
Zagier noted an one-to-one correspondence
$(\bZ+1/2) \times \cP \to \cX_0$
given by
$$(\beta, \lambda) \mapsto U(\beta) \ba_\lambda,$$
where $\cP$ denotes the set of all partitions.
This can be generalized to an one-to-one correspondence:
\be \label{eqn:Shift}
(\bZ+1/2+\alpha) \times \cP \to \cX_\alpha.
\ee

For those of whom do not feel comfortable working with infinite wedge products,
$\cF_\alpha$ can also be treated in the following way.
For $\ba \in \cX_\alpha$,
write $\ba_+=\{a_1, \dots, a_r\}$ and $\ba^c_-=\{-b_1, \dots, -b_s\}$,
where $b_1 < \cdots < c_s$ are nonnegative.
Define
\be
\hat{e}_{\ba} = e_{a_1} \wedge \cdots \wedge e_{a_r} \wedge e_{-b_1} \wedge \cdots \wedge e_{-b_s}.
\ee
Let $\hat{\cF}_\alpha = \Span \{\hat{e}_\ba\;|\; \ba \in \cX_\alpha\}$.
Then $\cF_{\alpha} \cong \hat{\cF}_\alpha$.

\subsection{Dirac flow and its asymptotic expansion}
We call the family of operators $\{e^{tD}\}_{t\in \bR}$ the {\em Dirac flow}.
\be
e^{tD} e_{n+\alpha}(x) = e^{(n+\alpha)t} \cdot e_{n+\alpha}(x).
\ee
We define the action of $e^{tD}$ on $\cF_\alpha$ as follows.
For $\ba \in \cX_\alpha$,
\ben
e^{tD} e_\ba = \sum_{i=1}^\infty e_{a_1} \wedge \cdots \wedge e^{tD}e_{a_i} \wedge \cdots.
\een
It is clear that
\ben
e^{tD} e_\ba = \sum_{i=1}^\infty e^{ta_i} \cdot e_\ba.
\een
For $\alpha \in [0, 1)$, we then have
\ben
W_\ba(t) & = & \sum_{i=1}^\infty e^{ta_i}
= \sum_{a_i \in \ba_+} e^{ta_i} - \sum_{a_i \in \ba^c_-} e^{ta_i}
+ \sum_{n=0}^\infty e^{-(n+\alpha) t} \\
& = & \sum_{a_i \in \ba_+} e^{ta_i} - \sum_{a_i \in \bb^c_-} e^{ta_i}
+ \frac{e^{-\alpha t}}{1 - e^{-t}} \\
& = & \sum_{n=0}^\infty \frac{t^n}{n!}
(\sum_{a_i \in \bb_+} a_i^n - \sum_{a_i \in \bb^c_-} a_i^n)
+ \frac{1}{t} + \sum_{k=1}^\infty \frac{B_{k}(1-\alpha)}{k!} t^{k-1}.
\een
So if one writes
\be \label{eqn:W}
W_\ba(t) = \sum_{k=0}^\infty Q_k(\ba) t^{k-1},
\ee
then
\bea
&& Q_0(\ba) = 1, \\
&& Q_k(\ba) = \frac{1}{(k-1)!} (\sum_{a_i \in \bb_+} a_i^{k-1}
- \sum_{a_i \in \bb^c_-} a_i^{k-1})
+ \frac{B_{k}(1-\alpha)}{k!}.
\eea
This gives us an asymptotic expansion of the Dirac flow on the Fock space $\cF_\alpha$
as $t\to 0$.
Recall the Hurwitz zeta function is defined by:
\be
\zeta(s, \theta) = \sum_{n=0}^\infty \frac{1}{(n + \theta)^s}.
\ee
And  its values at nonpositive integers are give by the following formula:
\be
\zeta(-n,\theta) = - \frac{B_{n+1}(\theta)}{n+1}.
\ee
So we have
\be
(k-1)!Q_k(\ba) = \sum_{a_i \in \ba_+} a_i^{k-1}
- \sum_{a_i \in \ba^c_-} a_i^{k-1}
- \zeta(-k+1, 1-\alpha).
\ee
We understand $Q_k(\ba)$ as the eigenvalue some operator $Q_k$ acting on $e_\ba$.
The operators $Q_k$ will be called the modes of the Dirac flow.

\begin{rmk}
Zagier \cite{Zag} considered $W_\ba$ and $Q_k(\ba)$ for $\alpha = 0$ and $\half$.
\end{rmk}

\subsection{Character of the modes of the Dirac flow}
Define
\be
\chi_\alpha(q_0, q_1, \dots) = \sum_{\ba \in \cX_{\alpha}} q_0^{Q_1(\ba)} q_1^{Q_2(\ba)} \cdots
= \tr e^{t_0Q_1+t_1Q_2+\cdots}|_{\cF_\alpha},
\ee
where $q_n = e^{t_n}$.
One can compute $\chi_\alpha$ in two different ways.
First of all,
using the isomorphism $\cF_\alpha \cong \hat{\cF}_\alpha$,
one can convert the computation to that of a trace of some operator
on the latter space,
which is easily evaluated to give us:
\be \label{eqn:Fermion}
\begin{split}
& \chi_\alpha(q_1, q_2, \dots) =
\exp(-\sum_{n=0}^\infty \frac{t_{n}}{(n-1)!} \zeta(-n, 1-\alpha)) \\
& \;\;\;\; \cdot \prod_{n \geq 1}
(1 + \exp(\sum_{j=0}^\infty \frac{t_j}{j!}(n-1+\alpha)^j))
(1 + \exp(- \sum_{j=0}^\infty \frac{t_j}{j!} (-n+\alpha)^j)).
\end{split}
\ee
Secondly, one can use the one-to-one correspondence in \eqref{eqn:Shift} as follows:
\ben
\chi_\alpha(q_1, q_2, \dots) = \sum_{n \in \bZ}
\sum_{\lambda \in \cP} q_1^{Q_1(\ba_\lambda+n+1/2+\alpha)}
q_2^{Q_2(\ba_\lambda+n+\half+\alpha)} \cdots
\een
One can use \eqref{eqn:W} to compute
$Q_k(a_\lambda+n+\half+\alpha)$ in the following way:
\ben
&& \sum_{k=0}^\infty Q_k(\ba_\lambda+n+\half+\alpha) t^{k-1} = W_{\ba_\lambda+n+\half+\alpha}(t)
= \sum_{i=1}^\infty e^{(\lambda_i-i+\half + n+ \half +\alpha)t} \\
& = & e^{(n+\half + \alpha)t} \cdot \sum_{i=1}^\infty e^{(\lambda_i-i+\half)t}
= e^{(n+\half + \alpha)t} \cdot W_{\ba_\lambda}(t) \\
& = & \sum_{l \geq 0} \frac{t^l}{l!}(n+\half+\alpha)^l \cdot \frac{1}{t}
\sum_{m=0}^\infty Q_m(\ba_\lambda) t^m,
\een
it follows that
\be
Q_k(\ba_\lambda+n+\half+\alpha)
= \sum_{l+m=k} \frac{1}{l!}(n+\half+\alpha)^l \cdot Q_m(\ba_\lambda).
\ee
Now we have:
\ben
&& \sum_{k=1}^\infty t_{k-1} Q_k(\ba_\lambda+n+\half+\alpha)
= \sum_{k=1}^\infty t_{k-1} \sum_{l+m=k} \frac{1}{l!}(n+\half+\alpha)^l \cdot Q_m(\ba_\lambda) \\
& = & \sum_{k=1}^\infty \frac{t_{k-1}}{k!}(n+\half+\alpha)^k
+ \sum_{m=1}^\infty Q_m(\ba_\lambda) \sum_{l=0}^\infty \frac{t_{l+m-1}}{l!}(n+\half+\alpha)^l.
\een
Therefore,
we get another formula for $\chi_\alpha$:
\be \label{eqn:Boson}
\begin{split}
& \chi_\alpha(q_1, q_2, \dots) \\
= &\sum_{n \in \bZ}
\sum_{\lambda \in \cP}
\exp\biggl( \sum_{k=1}^\infty \frac{t_{k-1}}{k!}(n+\half+\alpha)^k
+ \sum_{m=1}^\infty Q_m(\ba_\lambda) \sum_{l=0}^\infty \frac{t_{l+m-1}}{l!}(n+\half+\alpha)^l\biggr).
\end{split}
\ee
Now by comparing \eqref{eqn:Fermion} with \eqref{eqn:Boson},
we get the following formula:
\be \label{eqn:Boson-Fermion}
\begin{split}
&\sum_{n \in \bZ}
\sum_{\lambda \in \cP}
\exp\biggl( \sum_{k=1}^\infty \frac{t_{k-1}}{k!}(n+\half+\alpha)^k
+ \sum_{m=1}^\infty Q_m(\ba_\lambda) \sum_{l=0}^\infty \frac{t_{l+m-1}}{l!}(n+\half+\alpha)^l\biggr) \\
= & \exp(-\sum_{n=0}^\infty \frac{t_{n}}{n!} \zeta(-n, 1-\alpha)) \\
&  \cdot \prod_{n \geq 1}
(1 + \exp(\sum_{j=0}^\infty \frac{t_j}{j!}(n-1+\alpha)^j))
(1 + \exp(- \sum_{j=0}^\infty \frac{t_j}{j!} (-n+\alpha)^j)).
\end{split}
\ee

Let us look at some special case of this formula.
First take $\alpha = \half$ to get:
\be
\begin{split}
&\sum_{n \in \bZ}
\sum_{\lambda \in \cP}
\exp\biggl( \sum_{k=1}^\infty \frac{t_{k-1}}{k!}n^k
+ \sum_{m=1}^\infty Q_m(\ba_\lambda) \sum_{l=0}^\infty \frac{t_{l+m-1}}{l!}n^l\biggr) \\
= & \exp(-\sum_{n=0}^\infty \frac{t_{n}}{n!} \zeta(-n, \half)) \\
&  \cdot \prod_{n \geq 1}
(1 + \exp(\sum_{j=0}^\infty \frac{t_j}{j!}(n-\half)^j))
(1 + \exp(- \sum_{j=0}^\infty \frac{t_j}{j!} (-n+\half)^j)).
\end{split}
\ee
Next we take $t_2 = t_3 = \cdots = 0$,
and use the following identities \cite[(3)]{Zag}
\be
Q_1(\ba_\lambda) = 0, \;\;\; Q_2(\ba_\lambda) = |\lambda|- \frac{1}{24}
\ee
to get:
\be
\begin{split}
&\sum_{n \in \bZ}
\sum_{\lambda \in \cP}
\exp\biggl(nt_0 + \half n^2 t_1
+ (|\lambda|- \frac{1}{24})t_1 \biggr) \\
= & q_1^{-1/24}  \cdot \prod_{n \geq 1}
(1 + q_0 q_1^{n-\half})
(1 + q_0^{-1} q_1^{n-\half}).
\end{split}
\ee
After using the Euler identity:
\be
\sum_{\lambda \in \cP} q_1^{|\lambda|} = \frac{1}{\prod_{n=1}^\infty (1-q_1^n)},
\ee
we recover the Jacobi triple product identity:
\be
\sum_{n \in \bZ} q_0^nq_1^{\half n^2}
= \prod_{n\geq 1} (1-q_1^n)(1+ q_0 q_1^{n-\half})
(1 + q_0^{-1} q_1^{n-\half}).
\ee

Next, we take $\alpha = 0$ in \eqref{eqn:Boson-Fermion} to get:
\be
\begin{split}
&\sum_{n \in \bZ}
\sum_{\lambda \in \cP}
\exp\biggl( \sum_{k=1}^\infty \frac{t_{k-1}}{k!}(n+\half)^k
+ \sum_{m=1}^\infty Q_m(\ba_\lambda) \sum_{l=0}^\infty \frac{t_{l+m-1}}{l!}(n+\half)^l\biggr) \\
= & \exp(-\sum_{n=0}^\infty \frac{t_{n}}{n!} \zeta(-n))
\cdot (1+ \exp t_0) \\
&  \cdot \prod_{n \geq 1}
(1 + \exp(\sum_{j=0}^\infty \frac{t_j}{j!}n^j))
(1 + \exp(- \sum_{j=0}^\infty \frac{t_j}{j!} (-n)^j)).
\end{split}
\ee
After taking $t_0 = \pi \sqrt{-1}$,
we get:
\be
\begin{split}
& \sum_{n \in \bZ} (-1)^n
\sum_{\lambda \in \cP}
\exp\biggl( \sum_{k=2}^\infty \frac{t_{k-1}}{k!}(n+\half)^k
+ \sum_{m=2}^\infty Q_m(\ba_\lambda) \sum_{l=0}^\infty \frac{t_{l+m-1}}{l!}(n+\half)^l\biggr) \\
= & 0.
\end{split}
\ee
It can be rewritten in a more elegant form:
\be
\sum_{n \in \bZ} (-1)^n
\sum_{\lambda \in \cP}
\exp\biggl( \sum_{k=2}^\infty t_{k-1}
\sum_{l=0}^k  \frac{1}{l!} (n+\half)^lQ_{k-l}(\ba_\lambda)\biggr) \\
= 0.
\ee
It can be further simplified if one uses Zagier's operator ${\mathbf \pd}$ defined by:
\be
{\bf \pd} = \sum_{k=1}^\infty Q_{k-1}\frac{\pd}{\pd Q_k}.
\ee
Then the above recursion relation can be rewritten as
\be
\sum_{n \in \bZ} (-1)^n
\sum_{\lambda \in \cP}
\exp\biggl[\biggl(e^{(n+\half)\pd} \sum_{k=2}^\infty t_{k-1} Q_k \biggr)(\ba_\lambda)\biggr]
= 0.
\ee
This is just Zagier's recursion relations \cite[Theorem 1]{Zag} in the form
of a generating series.
Therefore,
\eqref{eqn:Boson-Fermion} is a generalization of both Jacobi triple product identity
and Zagier's recursion relations for $q$-brackets of shifted symmetric polynomials.

\end{appendix}

\vspace{.2in}

{\em Acknowledgements}.
The work in this note was started and completed during the author's attendance of a conference
on quantum K-theory and the stay after that at Sun Yat-Sen University.
The author thanks Professor Jianxun Hu and Professor Changzheng Li and other colleagues there
for the invitation and the hospitality. He also thanks
Professor Weiping Li, Professor Siqi Liu and Professor Yongbin Ruan for
stimulating discussions.
During the research in this work the author is partly supported by NSFC grant 11661131005.

\end{document}